\documentclass[12pt,leqno]{amsart}
\makeatletter
\usepackage{amsfonts,delarray,amssymb,amsmath,amsthm,a4,a4wide}
\usepackage[top=2.5cm, bottom = 2.5cm, left = 2.5cm, right=2.5cm]{geometry}
\usepackage{epsfig}
\usepackage{color}

\title[On the second inner variation of the Allen-Cahn Functional]{On the second inner variation of the Allen-Cahn Functional and its applications
}
\author{Nam Q. Le}
\address{Department of
Mathematics, Columbia University, New York,
 USA}
\email{namle@math.columbia.edu}
\newcommand{\review}[2][\right]{\relax
\ifx#1\right\relax \left.\fi#2#1\rvert}
\let\abs=\envert
 
\newtheorem{theorem}{Theorem}
\newtheorem{claim}{Claim}[section]

\newtheorem{remark}{Remark}[section]

\newcommand{\bef}{\begin{flushright}}
\newcommand{\eef}{\end{flushright}}

\newcommand{\eval}[2][\right]{\relax
\ifx#1\right\relax \left.\fi#2#1\rvert}

\let\abs=\envert

\numberwithin{equation}{section}

\newcommand\e{\varepsilon}
\newcommand{\h}{\hspace*{.24in}}

\def\h{\hspace*{.24in}}

\def\beq{\begin{eqnarray*}}
\def\eeq{\end{eqnarray*}}

\def\RR{\mbox{$I\hspace{-.06in}R$}}

\newenvironment{myindentpar}[1]%
{\begin{list}{}%
         {\setlength{\leftmargin}{#1}}%
         \item[]%
}
{\end{list}}

\begin{document}
\date{September 28, 2010}
\maketitle

\pagenumbering{arabic}

\begin{abstract}
In this paper, we study the relation between the second inner variations of the Allen-Cahn functional and its Gamma-limit, the area functional. Our result
implies that the Allen-Cahn functional only approximates well the area functional up to the first order. 
However, as an application of our result, we prove, assuming the single-multiplicity property of the limiting energy, that
the Morse indices of critical points of the Allen-Cahn functional are bounded from below by the Morse index of the limiting minimal hypersurface.
\end{abstract}
\section{Introduction and Main Results}
\h Let $\Omega$ be an open smooth bounded set in $\RR^{N}$ ($N\geq 2$). Then, for any $C^{2}$, closed hypersurface $\Gamma$ inside $\Omega$ with finite perimeter, we can use the Allen-Cahn functional to approximate the area of $\Gamma$. Indeed, for each $\varepsilon >0$, consider the following Allen-Cahn functional
\begin{equation} E_{\varepsilon}(u)=\int_{\Omega}\frac{\varepsilon \abs{\nabla u}^2}{2} +\frac{W(u)}{\varepsilon},\label{MM}\end{equation} 
where $W(u) = \frac{1}{2} (1-u^2)^2$ is the double-well potential and $u: \RR^{N}\rightarrow \RR$ is a scalar function. This is a typical 
energy modeling the phase separation phenomena within the van der Waals-Cahn-Hilliard gradient theory of phase transitions \cite{AC}. Then, we can find a sequence of scalar 
functions $u^{\varepsilon}$ such that 
\begin{myindentpar}{1cm}

\begin{equation}\text{The zero level sets of~} u^{\e} \text{converge to~} \Gamma \text{~in the Hausdorff distance sense}
 \label{zeroH}
\end{equation}

and 
\begin{equation}
 \lim_{\varepsilon \rightarrow 0} E_{\varepsilon}(u^{\varepsilon}) = 2\sigma \mathcal{H}^{N-1}(\Gamma),
\label{appro}
\end{equation}
\end{myindentpar}

where $\sigma =\displaystyle \int_{-1}^{1}\sqrt{W(s)/2} ds=\frac{2}{3}$ and $\mathcal{H}^{N-1}$ denotes the $(N-1)$-dimensional Hausdorff measure. There are 
many such sequences of $u^{\varepsilon}$; the construction of one such sequence follows from the construction part of the general result in 
the framework of $\Gamma$-convergence (see Modica-Mortola \cite{MM} and Sternberg \cite{Stern}, for example). Note also that $E_{\varepsilon}$ $\Gamma$-converges to the area functional 
\begin{equation}E(u)=2\sigma \mathcal{H}^{N-1}(\Gamma) := E(\Gamma).\end{equation} 
\h Here $u$ is a function of bounded variation taking values $\pm 1$, $\Gamma$ is the interface separating 
the phases, i.e, $\Gamma = \partial\{x\in \Omega: u(x)=1\}\cap \Omega$. Roughly speaking, $\Gamma$ is the limit of the zero level sets of $u^{\e}$. 
Furthermore, we know 
that (\ref{zeroH}) and (\ref{appro}) imply the single-multiplicity property, i.e., in the sense of 
Radon measures, \begin{equation}\left(\frac{\varepsilon\abs{\nabla u^{\varepsilon}}^2}{2} +\frac{W(u^{\varepsilon})}{\varepsilon}\right)dx 
\rightharpoonup 2\sigma d\mathcal{H}^{N-1}\lfloor\Gamma. \label{single} \end{equation} 
Therefore, from the work of Reshetnyak \cite{Res}, we can prove that the first inner variation of $E_{\varepsilon}$ at $u^{\varepsilon}$
\begin{equation}
 \delta E_{\varepsilon}(u^{\varepsilon},\eta) = \int_{\Omega}\left( \frac{\varepsilon \abs{\nabla u^{\varepsilon}}^2}{2} +\frac{W(u^{\varepsilon})}{\varepsilon}\right)\text{div} \eta - \varepsilon(\nabla u^{\varepsilon},\nabla u^{\varepsilon}\cdot\nabla\eta),
\label{fvep}\end{equation}
converges to the first inner variation of $E$ at $\Gamma$
\begin{equation}
 \delta E(\Gamma,\eta) = 2\sigma\int_{\Gamma} (\text{div}\eta -
\partial_{k}\varphi^{j}n_{j}\cdot n_{k})d\mathcal{H}^{N-1}.
\label{fve}\end{equation}
In (\ref{fvep}) and (\ref{fve}), $\eta\in (C_{c}^{1}(\Omega))^{N}$ is a vector field and 
$\stackrel{\rightarrow}{n} = (n_{1},\cdots,n_{N})$ denotes the outward unit normal to the region enclosed by $\Gamma$; and $(\cdot,\cdot)$ denotes the standard inner product on $R^{N}$.\\ 
\h This convergence result especially imposes the criticality conditions on $\Gamma$, typically stationary or minimal condition, if one would like 
to approximate the area of $\Gamma$ by $u^{\varepsilon}$ which are critical points of $E_{\varepsilon}$. 
The most general result concerning the geometric properties of $\Gamma$ as the limit of the zero level sets of $u^{\e}$ with
 suitable uniform Sobolev bounds on $\e \Delta u^{\e} -\e^{-1} W^{'} (u^{\e})$ is due to Tonegawa
\cite{ToneDiff}. Regarding $L^{2}$-bounds on $\e \Delta u^{\e} -\e^{-1} W^{'} (u^{\e})$, one could mention the work of R\"{o}ger and Sch\"{a}tzle on De Giorgi's conjecture. See \cite{RS}
and the references therein.\\
\h If $u^{\varepsilon}$ (resp. $\Gamma$) are critical points of $E_{\varepsilon}$ (resp. $E$), then a natural question to ask is: what are the relations 
between the stability of $u^{\varepsilon}$ with respect to $E_{\varepsilon}$ and that of $\Gamma$ with respect to $E$? 
Assuming the stability of $u^{\e}$, by using clever test functions in the stability inequality satisfied by $u^{\e}$, Tonegawa \cite{Tone} proved
that $\Gamma$ must be a stable varifold. Here no assumptions on the regularity of $\Gamma$ are assumed a priori. 
Another way to answer the above question is to study the second variations of both functionals $E_{\varepsilon}$ and $E$. We find that, contrary 
to the first variation, the second variation of $E_{\varepsilon}$ at $u^{\varepsilon}$, $\delta^{2} E_{\varepsilon}(u^{\varepsilon},\eta,\zeta)$, does not in general 
converges to the second variation of $E$ at $\Gamma$, $\delta^{2}E(\Gamma,\eta,\zeta)$! \\
\h The purpose of this note is two-fold. First, under some regularity assumptions, we provide the precise relation between the limit of the second inner variation of $E_{\e}$, assuming the single-multiplicity
condition (\ref{appro}), and the second inner variation of $E$. This relation can be of independent interest. Second, we use this relation to estimate from below
the Morse indices of the critical points
$u^{\e}$ of $E_{\e}$ for $\e$ sufficiently small in terms of the Morse index of the critical point $\Gamma$ of $E$. Here, again, $\Gamma$ is the limit of the zero
level set of $u^{\e}$.\\
\h Before stating our main result, we recall some standard definitions. Consider smooth vector fields $\eta,\zeta\in (C_{c}^{1}(\Omega))^{N}$. Then, for $t$ sufficiently small, the map $\Phi_{t} (x) = x + t\eta(x) + \frac{t^2}{2}\zeta(x)$ is a diffeomorphism of $\Omega$ into itself. We think of $\eta$ and $\zeta$ as initial velocity and acceleration vectors when we deform the domain $\Omega$. The second inner variation of $E_{\e}$ at $u^{\e}$ with respect to the velocity and acceleration vectors $\eta$ and $\zeta$ is defined by
\begin{equation*}
 \delta^{2} E_{\varepsilon}(u^{\varepsilon},\eta,\zeta) = \left.\frac{d^2}{dt^2}\right\rvert_{t=0} E_{\varepsilon}(u^{\varepsilon}_{t}),
\end{equation*}
where $u^{\varepsilon}_{t}(y) = u^{\varepsilon}(\Phi^{-1}_{t}(y))$. The second inner variation of $E$ at $\Gamma$ with respect to the velocity and acceleration vectors $\eta$ and $\zeta$ is defined by
\begin{equation*}
 \delta^{2}E(\Gamma,\eta,\zeta) = \left.\frac{d^2}{dt^2}\right\rvert_{t=0} E(\Gamma_{t}),
\end{equation*}
where $\Gamma_{t} = \Phi_{t}(\Gamma)$.\\
\h In this note, we prove the following main result, revealing the exact discrepancy between the second variation of $E_{\varepsilon}$ at $u^{\varepsilon}$ and that of $E$ at $\Gamma$.
\begin{theorem}
 Let $\Gamma$ be any $C^{2}$, closed hypersurface inside $\Omega$ with finite perimeter and let $u^{\varepsilon}$ be any sequence of scalar functions such that
the zero level sets of $u^{\e}$ converge to $\Gamma$ in the Hausdorff distance sense and 
\begin{equation}
 \lim_{\varepsilon \rightarrow 0} E_{\varepsilon}(u^{\varepsilon}) = 2\sigma \mathcal{H}^{N-1}(\Gamma).
\label{appro2}
\end{equation} Then, for all smooth vector fields $\eta,\zeta\in (C_{c}^{1}(\Omega))^{N}$, we 
have \begin{equation}\lim_{\varepsilon\rightarrow 0}\delta^{2} E_{\varepsilon}(u^{\varepsilon},\eta,\zeta) =
 \delta^{2}E(\Gamma,\eta,\zeta) + 2\sigma\int_{\Gamma} (\stackrel{\rightarrow}{n},\stackrel{\rightarrow}{n}\cdot\nabla\eta)^2.\label{discrep}\end{equation} 
\label{mainthm}
\end{theorem}
\begin{remark}
The discrepancy term is $2\sigma\int_{\Gamma} (\stackrel{\rightarrow}{n},\stackrel{\rightarrow}{n}\cdot\nabla\eta)^2.$ It is new and has a sign.
Our result is proved without assuming any criticality conditions on $u^{\varepsilon}$ nor $\Gamma$. Thus, on the levels of energy and the first inner variation, the Allen-Cahn functionals approximate well the area functional. This is no longer true for the second inner variation.
\end{remark}
\begin{remark}
 The formula for $\delta^{2}E(\Gamma,\eta,\zeta)$ is given by (\ref{sve}). Let $\zeta\equiv 0$ and $\eta$ be a normal vector field defined on $\Gamma$, i.e.,
$\eta = f \stackrel{\rightarrow}{n}$ for some function $f\in C_{c}^{1}(\Omega)$. Then (\ref{discrep}) and (\ref{sve}) give
\begin{equation}
 \lim_{\varepsilon\rightarrow 0}\delta^{2} E_{\varepsilon}(u^{\varepsilon},f \stackrel{\rightarrow}{n},0) = \int_{\Gamma} \abs{\nabla f}^2 - \abs{A}^2 f^2
\label{fullgrad}
\end{equation}
where $\abs{A}$ denotes the length of the second fundamental form $A$ of $\Gamma$. The quantity on the right hand side of (\ref{fullgrad}) is the one used by
Tonegawa in his stability result, \cite[Theorem 3]{Tone}. For the proof of the stability of the interface $\Gamma$, this is sufficient. Our formula 
(\ref{fullgrad}) explains that in general, we cannot replace the full gradient $\int_{\Gamma} \abs{\nabla f}^2$ by the restricted gradient 
$\int_{\Gamma} \abs{\nabla^{\Gamma} f}^2$ in Tonegawa's stability result. 
\end{remark}

Let us denote by $D^{2}E_{\e}(u)$ the Hessian of $E_{\e}$ at $u$ and $Q_{\e}(u)$ the associated quadratic function, associated to the bilinear
continuous function $B_{\e}(u)(\cdot,\cdot)$. If $u_{t}$ is a variation of $u$, i.e, $u_{0}= u,$ then
\begin{equation*}
 \left.\frac{d^2}{dt^2}\right\rvert_{t=0} E_{\varepsilon}(u_{t}) = Q_{\e}(u)( \left.\frac{d}{dt}\right\rvert_{t=0} u_{t}) = B_{\e}(u)
(\left.\frac{d}{dt}\right\rvert_{t=0} u_{t}, \left.\frac{d}{dt}\right\rvert_{t=0}u_{t} ). 
\end{equation*}
Similarly, we can define $D^{2} E(u), Q(u)$ and $B(u)$ for $E$. 
Now, let $\zeta\equiv 0$ and let $\eta $ be a normal vector field defined on $\Gamma$. 
Assuming the smoothness of $\Gamma$, we can find an extension $\tilde{\eta}$ of $\eta$ to $\Omega$ such that $(\stackrel{\rightarrow}{n}, 
\stackrel{\rightarrow}{n}\cdot\nabla\tilde{\eta}) =0$. In this case, combining 
(\ref{discrep}) with Theorem 1.1 in \cite{Serfaty}, we obtain the following result.
\begin{theorem}
Let $u^{\varepsilon}$ be critical 
points of $E_{\varepsilon}$ , i.e., $\delta E_{\varepsilon}(u^{\varepsilon},\eta) = 0$. Upon extracting a subsequence, $u^{\e}\longrightarrow u \in BV(\Omega,
\{1,-1\})$ in $L^{1}(\Omega)$. Let $\Gamma$ be the interface separating the phases of $u$. Assume that (\ref{appro2}) holds and $\Gamma$ is of class
$C^{2}$. Denote by $n_{\e}^{+}$ 
the dimension (possibly infinite) of the space spanned by the eigenvectors of $D^{2}E_{\e}(u_{\e})$ associated to positive eigenvalues, and $n^{+}$ 
the dimension (possibly infinite) of the space spanned by the eigenvectors of $D^{2}E(u)$ associated to positive eigenvalues
(resp. $n^{-}_{\e}$ and $n^{-}$ for negative eigenvalues). Then, for $\e$ small enough we have 
\begin{equation}
 n^{+}_{\e}\geq n^{+},~n^{-}_{\e}\geq n^{-}.
\label{MI}
\end{equation}
\label{Morse}
\end{theorem}
We note that by writing $\tilde{\eta} = f \stackrel{\rightarrow}{n}$, we have
\begin{equation*}
 Q(u)(\tilde{\eta}) = \int_{\Gamma} \abs{\nabla^{\Gamma} f}^2 - \abs{A}^2 f^2
\end{equation*}
where $\abs{A}$ denotes the length of the second fundamental form $A$ of $\Gamma$. Therefore,
Theorem \ref{Morse} implies, in particular, that if $n^{-}_{\e}$ is the Morse index of $u_{\e}$ then $\Gamma$ 
must be a generalized minimal hypersurface with Morse index $n^{-}$
satisfying
\begin{equation}
 \liminf_{\e\rightarrow 0}n^{-}_{\e}\geq n^{-}.
\label{MorseIneq}
\end{equation}
Thus, if $u_{\e}$ is stable then $\Gamma$ is stable, reproving a special case of Tonegawa's stability result when
$\Gamma$ is smooth and (\ref{appro2}) is satisfied. Regarding stability theory, see Tonegawa \cite{Tone} for a very general result without assuming (\ref{appro2})
; see also Serfaty \cite{Serfaty} for 
the complex-valued version of $u^{\varepsilon}$. Regarding regularity theory for stable hypersurface $\Gamma$, see a very interesting recent paper by 
Tonegawa and Wickramasekera \cite{TW}. 
\begin{remark}
 When $\Gamma$ is a minimal hypersurface satisfying certain nondegeneracy conditions, Pacard and Ritor\'{e} \cite{PR} constructed critical points $u^{\e}$ of $E_{\e}$ whose zero 
level sets converge to $\Gamma$ and (\ref{appro2}) holds. Thus, Theorem \ref{Morse} provides an estimate for the Morse indices of $u^{\e}$ constructed 
in \cite{PR} in terms of the Morse
index of $\Gamma$. 
\end{remark}
\begin{remark}
If $N=3$, $\Omega =\RR^{3}$ and $\Gamma$ is a minimal surface, embedded, complete with finite total curvature and is nondegenerate, then del Pino, Kowalczyk and 
Wei \cite{dPKW} constructed critical points $u^{\e}$ of $E_{\e}$ whose zero level sets converge to $\Gamma$ and $n_{\e}^{-}(u^{\e}) = n^{-}(\Gamma)$ for $\e$
sufficiently small. 
\end{remark}
\begin{remark}
 In general, the inequality (\ref{MorseIneq}) can be strict. Here is one example for $N=2$ with $\Gamma$ singular at one point. Let
$\Omega$ be the unit disc in $\RR^{2}$ and $\Gamma$ be the cross 
$\Gamma = \{(x_{1}, 0), -1\leq x_{1}\leq 1\}\cup \{(0, x_{2}), -1\leq x_{2}\leq 1\}$. Using the construction of saddle solutions in Dang, Fife and 
Peletier \cite{DFP} (see also Gui \cite[Proposition 3.1]{Gui}), one can construct a sequence of critical points $u^{\e}$ of $E_{\e}$ such that 
the limit of the zero level set of $u^{\e}$
is the cross $\Gamma$. The cross is a stable varifold and thus has Morse index $0$. On the other hand, for $\e$ sufficiently small, the critical points $u^{\e}$
of $E_{\e}$ are not stable, and thus have Morse index at least 1. The reason that $u^{\e}$ are not stable is as follows. If otherwise, then by 
Tonegawa's result \cite[Theorem 5]{Tone}, the limit zero level set $\Gamma$ is a finite number of lines with no intersections or junctions. Therefore, it cannot be the cross, which 
is a contradiction.
\end{remark}
It would be very interesting to provide estimates similar to (\ref{MI}) when the multiplicity one condition (\ref{appro2}) is dropped. In this regard,
we have the following partial result, where we replace (\ref{appro2}) by the following mild conditions:
\begin{myindentpar}{1cm}
 (C1) The limit measure of $\left(\frac{\varepsilon\abs{\nabla u^{\varepsilon}}^2}{2} +\frac{W(u^{\varepsilon})}{\varepsilon}\right)dx$ is concentrated
on $\Gamma$.\\
(C2) $\Gamma$ is connected.
\end{myindentpar}
\begin{theorem}
Let $u^{\varepsilon}$ be critical 
points of $E_{\varepsilon}$ , i.e., $\delta E_{\varepsilon}(u^{\varepsilon},\eta) = 0$. Upon extracting a subsequence, $u^{\e}\longrightarrow u \in BV(\Omega,
\{1,-1\})$ in $L^{1}(\Omega)$. Let $\Gamma$ be the interface separating the phases of $u$. Assume that (C1) and (C2) are satisfied and $\Gamma$ is of class
$C^{2}$. Denote by $n_{\e}^{+}$ 
the dimension (possibly infinite) of the space spanned by the eigenvectors of $D^{2}E_{\e}(u_{\e})$ associated to positive eigenvalues, and $n^{+}$ 
the dimension (possibly infinite) of the space spanned by the eigenvectors of $D^{2}E(u)$ associated to positive eigenvalues
(resp. $n^{-}_{\e}$ and $n^{-}$ for negative eigenvalues). Then, for $\e$ small enough we have 
\begin{equation}
 n^{+}_{\e}\geq n^{+},~n^{-}_{\e}\geq n^{-}.
\label{mMI}
\end{equation}
\label{mMorse}
\end{theorem}

{\bf Acknowledgements:} The author would like to thank Yoshihiro Tonegawa for useful discussions and interesting suggestions during the preparation of this 
note.
\section{Proof of the Main Results}
\h This section is entirely devoted to the proof of Theorems \ref{mainthm}, \ref{Morse} and \ref{mMorse}. 
\begin{proof}[Proof of Theorem \ref{mainthm}]
First of all, we have the following formula for the second inner variation of $E$ at $\Gamma$ (see Simon \cite[p. 51]{Simon}, for example)
\begin{equation}
 \delta^{2}E(\Gamma,\eta,\zeta) = 2\sigma\int_{\Gamma}\left\{ \text{div}^{\Gamma}\zeta + (\text{div}^{\Gamma}\eta)^2 + \sum_{i=1}^{N-1}\abs{(D_{\tau_{i}}\eta)^{\perp}}^2 - \sum_{i,j=1}^{N-1}(\tau_{i}\cdot D_{\tau_{j}}\eta)(\tau_{j}\cdot D_{\tau_{i}}\eta)\right\},
\label{sve}\end{equation}
where $\text{div}^{\Gamma}\varphi$ denotes the tangential divergence of $\varphi$ on $\Gamma$; and for each point $x\in\Gamma$, $\{\tau_{1}(x),\cdots,\tau_{N-1}(x)\}$ is any orthonormal basis for the tangent space $T_{x}(\Gamma)$; for each $\tau\in T_{x}(\Gamma)$, $D_{\tau} \eta$ is the directional derivative and the normal part of $D_{\tau_{i}} \eta$ is denoted by \begin{equation*}
 (D_{\tau_{i}} \eta)^{\perp} = D_{\tau_{i}} \eta -\sum_{j=1}^{N-1}(\tau_{j}\cdot D_{\tau_{i}} \eta)\tau_{j}. 
\end{equation*}
Next, for the second inner variation of $E_{\varepsilon}$ at $u^{\varepsilon}$, we claim that
\begin{multline}
\delta^{2} E_{\varepsilon}(u^{\varepsilon},\eta,\zeta) =\int_{\Omega} \left\{ \left( \frac{\varepsilon \abs{\nabla u^{\varepsilon}}^2}{2} +
\frac{W(u^{\varepsilon})}{\varepsilon}\right) \left(\text{div}\zeta + (\text{div}\eta)^2 -
\text{trace}((\nabla \eta)^2)\right) \right.\\ + \varepsilon \abs{\nabla u^{\varepsilon}\cdot\nabla\eta}^2 + 
2\varepsilon (\nabla u^{\varepsilon},\nabla u^{\varepsilon}\cdot (\nabla\eta)^2) -
\e(\nabla u^{\e},\nabla u^{\e}\cdot\nabla\zeta)\\ \left.-2\varepsilon (\nabla u^{\varepsilon},\nabla u^{\varepsilon}\cdot\nabla\eta)\text{div}\eta\right\}.
\label{svep}\end{multline}
\h We indicate how to derive this formula. Let $\eta\in (C_{c}^{1}(\Omega))^{N}$ and $\zeta\in (C_{c}^{1}(\Omega))^{N}$ be vector fields and $t\neq 0$ sufficiently small such that the map $\Phi_{t} (x) = x + t\eta(x) + \frac{t^2}{2}\zeta(x)$ is a diffeomorphism of $\Omega$ into itself. Set  $u^{\varepsilon}_{t}(y) = u^{\varepsilon}(\Phi^{-1}_{t}(y))$. We are going to calculate
\begin{equation}
 \delta^{2} E_{\varepsilon}(u^{\varepsilon},\eta,\zeta) = \left.\frac{d^2}{dt^2}\right\rvert_{t=0} E_{\varepsilon}(u^{\varepsilon}_{t}).
\end{equation}
By change of variables $y=\Phi_{t}(x)$, we have 
\begin{equation}
 E_{\varepsilon}(u^{\varepsilon}_{t}) = \int_{\Omega}\left[\frac{\varepsilon \abs{\nabla u^{\varepsilon}\cdot\nabla \Phi_{t}^{-1}(\Phi_{t}(x))}^2}{2} + \frac{(1-\abs{u^{\varepsilon}(x)}^{2})^2}{2\varepsilon}\right]\abs{\text{det}\nabla \Phi_{t}(x)}dx.
\label{vari}
\end{equation}
We need to expand the right-hand side of the above formula up to the second power of $t$. For this purpose, we use the following identity for matrices $A$ and $B$
\begin{equation}
 \text{det}(I + tA + \frac{t^{2}}{2}B) = 1 + t\text{trace}(A) + \frac{t^2}{2}[\text{trace}(B) + (\text{trace}(A))^2 - \text{trace}(A^2)] + O(t^3).
\end{equation}
Therefore, 
\begin{multline}
\text{det}\nabla\Phi_{t}(x) =\text{det} (I + t\nabla \eta (x) +\frac{t^2}{2}\nabla \zeta) \\= 1 +  t\text{div} \eta + \frac{t^2}{2}[ \text{div}\zeta + (\text{div}\eta)^2 - \text{trace}((\nabla\eta)^2)] + O(t^3).
\label{deter}
\end{multline}
Note that
\begin{equation}
 \nabla\Phi_{t}^{-1}(\Phi_{t}(x)) = [I + t\nabla\eta (x) +\frac{t^2}{2}\nabla\zeta(x)]^{-1} = I - t\nabla\eta -\frac{t^2}{2}\nabla\zeta(x)+ t^{2}(\nabla\eta)^2 + O(t^3).
\label{inverse}
\end{equation}
Plugging (\ref{deter}) and (\ref{inverse}) into (\ref{vari}), we get (\ref{svep}) after some simple calculations. For the sake of completeness, we include here the calculations. First, we note that, for $t$ sufficiently small, $\text{det}\nabla\Phi_{t}(x)>0$ and thus $\abs{\text{det}\nabla\Phi_{t}(x)} = \text{det}\nabla\Phi_{t}(x).$ Second, from (\ref{inverse}), we find that
\begin{equation*}
 \nabla u^{\varepsilon}\cdot\nabla \Phi_{t}^{-1}(\Phi_{t}(x)) = \nabla u^{\e} - t\nabla u^{\e}\cdot\nabla\eta -\frac{t^2}{2}\nabla u^{\e}\cdot\nabla\zeta(x)+ t^{2}\nabla u^{\e}\cdot(\nabla\eta)^2 + O(t^3).
\end{equation*}
Hence
\begin{multline}
 \frac{\varepsilon \abs{\nabla u^{\varepsilon}\cdot\nabla \Phi_{t}^{-1}(\Phi_{t}(x))}^2}{2} \\= \frac{\e}{2}\left\{ \abs{\nabla u^{\e}}^2 -2t (\nabla u^{\e}, \nabla u^{\e}\cdot\nabla \eta) + t^2 \abs{\nabla u^{\e}\cdot\nabla\eta}^2\right. \\+ \left. 2t^2(\nabla u^{\e}, \nabla u^{\e}\cdot (\nabla\eta)^2) -t^2 (\nabla u^{\e},\nabla u^{\e}\cdot\nabla\zeta) + O(t^3)\right\}.
\end{multline}
It follows that
\begin{multline}
 \left[\frac{\varepsilon \abs{\nabla u^{\varepsilon}\cdot\nabla \Phi_{t}^{-1}(\Phi_{t}(x))}^2}{2} + 
\frac{(1-\abs{u^{\varepsilon}(x)}^{2})^2}{2\varepsilon}\right]\abs{\text{det}\nabla \Phi_{t}(x)}\\ = 
\left(\frac{\e}{2}\left\{ \abs{\nabla u^{\e}}^2 -2t (\nabla u^{\e}, \nabla u^{\e}\cdot\nabla \eta) + 
t^2 \abs{\nabla u^{\e}\cdot\nabla\eta}^2 + 2t^2(\nabla u^{\e}, \nabla u^{\e}\cdot (\nabla\eta)^2) -t^2 (\nabla u^{\e},\nabla u^{\e}\cdot\nabla\zeta)\right\}\right.\\ +
 \left.\frac{(1-\abs{u^{\varepsilon}(x)}^{2})^2}{2\varepsilon} + O(t^3) \right)\left(1 +  t\text{div} \eta + \frac{t^2}{2}[ \text{div}\zeta + 
(\text{div}\eta)^2 - \text{trace}((\nabla\eta)^2)] + O(t^3)\right)\\= (\frac{\e}{2} \abs{\nabla u^{\e}}^2 + 
\frac{(1-\abs{u^{\varepsilon}(x)}^{2})^2}{2\varepsilon})\frac{t^2}{2}[ \text{div}\zeta + (\text{div}\eta)^2 - \text{trace}((\nabla\eta)^2)] \\ + 
\frac{\e t^2}{2}\{ \abs{\nabla u^{\e}\cdot\nabla\eta}^2 + 2(\nabla u^{\e}, \nabla u^{\e}\cdot (\nabla\eta)^2) -
(\nabla u^{\e},\nabla u^{\e}\cdot\nabla\zeta)\} \\ - \e t^2 (\nabla u^{\e}, \nabla u^{\e}\cdot\nabla \eta) \text{div}\eta\\ + \text{lower order terms in $t$} + O(t^3).
\end{multline}
Substituting this relation into (\ref{vari}), one gets (\ref{svep}) as desired.\\
 \h From (\ref{appro2}) and the fact that the zero level sets of $u^{\e}$ converge to $\Gamma$ in the Hausdorff distance sense, we have the single-multiplicity (\ref{single}). Consequently, from the work of Reshetnyak \cite{Res}, we can prove that
\begin{equation}
 \varepsilon \nabla u^{\varepsilon}\otimes\nabla u^{\varepsilon} dx\rightharpoonup 2\sigma \stackrel{\rightarrow}{n}\otimes\stackrel{\rightarrow}{n} 
\mathcal{H}^{N-1}\lfloor\Gamma.\label{res}
\end{equation}
(For a simple proof of this result, see Luckhaus and Modica \cite{LM}).\\ 
\h Passing to the limit in (\ref{svep}), employing (\ref{appro2}) and (\ref{res}), we obtain
\begin{multline}
 \lim_{\varepsilon\rightarrow 0}\delta^{2}E_{\varepsilon}(u^{\varepsilon},\eta,\zeta) = 2\sigma\int_{\Gamma}\text{div}\zeta + 
(\text{div}\eta)^2 - \text{trace} ((\nabla \eta)^2) \\+ 2\sigma\int_{\Gamma}\abs{\stackrel{\rightarrow}{n}\cdot\nabla\eta}^2 + 
2(\stackrel{\rightarrow}{n},\stackrel{\rightarrow}{n}\cdot(\nabla\eta)^2) -(\stackrel{\rightarrow}{n}, \stackrel{\rightarrow}{n}\cdot\nabla\zeta) - 
2 (\stackrel{\rightarrow}{n},\stackrel{\rightarrow}{n}\cdot\nabla\eta)\text{div}\eta.
\end{multline}
Note that $
\text{div}^{\Gamma}\eta = \text{div}\eta - (\stackrel{\rightarrow}{n},\stackrel{\rightarrow}{n}\cdot\nabla\eta). $
Hence
                                                                                                          \begin{multline}
\lim_{\varepsilon\rightarrow 0}\delta^{2}E_{\varepsilon}(u^{\varepsilon},\eta,\zeta) = 
2\sigma\int_{\Omega} \text{div}^{\Gamma}\zeta + (\text{div}^{\Gamma}\eta)^2 - \text{trace} ((\nabla \eta)^2) \\+ 
2\sigma\int_{\Gamma}\abs{\stackrel{\rightarrow}{n}\cdot\nabla\eta}^2 + 2(\stackrel{\rightarrow}{n},\stackrel{\rightarrow}{n}\cdot(\nabla\eta)^2) - 
(\stackrel{\rightarrow}{n},\stackrel{\rightarrow}{n}\cdot\nabla\eta)^2.
\label{svreduced}
\end{multline}                                                                     
Some calculation using local coordinates completes the proof of (\ref{discrep}). For the reader's convenience, we include the details. We can 
choose local coordinates so that $\{\tau_{1},\cdot,\tau_{N-1}, \stackrel{\rightarrow}{n}\}$ is the orthonormal basis of $R^{N}$. 
Furthermore, $\stackrel{\rightarrow}{n} = (0,\cdots, 0,1)$. We calculate successively\\
(i) $(\nabla\eta)_{ij} = \frac{\partial\eta^{i}}{\partial x_{j}}$,\\
(ii) $((\nabla\eta)^2)_{ij} = \sum_{k}\frac{\partial\eta^{i}}{\partial x_{k}}\frac{\partial\eta^{k}}{\partial x_{j}}$,\\
(iii) $\text{trace} (\nabla\eta)^2 = \sum_{i} ((\nabla\eta)^2)_{ii} = \sum_{i,k}\frac{\partial\eta^{i}}{\partial x_{k}}\frac{\partial\eta^{k}}{\partial x_{i}}$,\\
(iv) $2 (\stackrel{\rightarrow}{n}, \stackrel{\rightarrow}{n}\cdot (\nabla\eta)^2) = 2\sum_{i,j}n_{i}n_{j}((\nabla\eta)^2)_{ij} = 2 ((\nabla\eta)^2)_{NN} =2\sum_{k}\frac{\partial\eta^{N}}{\partial x_{k}}\frac{\partial\eta^{k}}{\partial x_{N}} $,\\
(v)$(\stackrel{\rightarrow}{n},\stackrel{\rightarrow}{n}\cdot\nabla\eta)^2 = (\sum_{i,j}n_{i}n_{j}\frac{\partial\eta^{i}}{\partial x_{j}})^2 = (\frac{\partial\eta^{N}}{\partial x_{N}})^2$,\\
(vi) $\abs{\stackrel{\rightarrow}{n}\cdot\nabla\eta}^2 = \abs{(\frac{\partial\eta^{j}}{\partial x_{i}}n_{j})}^2 = \sum_{i} \abs{\frac{\partial\eta^{N}}{\partial x_{i}}}^2$,\\
(vii) $\abs{\stackrel{\rightarrow}{n}\cdot\nabla\eta}^2 -(\stackrel{\rightarrow}{n},\stackrel{\rightarrow}{n}\cdot\nabla\eta)^2 = \sum_{i<N} \abs{\frac{\partial\eta^{N}}{\partial x_{i}}}^2$,\\
(viii) $(D_{\tau_{i}}\eta)^{\perp} = D_{\tau_{i}} \eta -\sum_{j=1}^{N-1}(\tau_{j}\cdot D_{\tau_{i}} \eta)\tau_{j} = (\frac{\partial\eta^{1}}{\partial x_{i}}, \cdots,\frac{\partial\eta^{N}}{\partial x_{i}}) -\sum_{j<N}\frac{\partial\eta^{j}}{\partial x_{i}}\tau_{j} = (0,\cdots,0, \frac{\partial\eta^{N}}{\partial x_{i}})$,\\
(ix) $\sum_{i<N}\abs{(D_{\tau_{i}}\eta)^{\perp}}^2 = \sum_{i<N}\abs{\frac{\partial\eta^{N}}{\partial x_{i}}}^2,$\\
(x) $\tau_{i}\cdot D_{\tau_{j}}\eta = \frac{\partial\eta^{i}}{\partial x_{j}}$,\\
(xi) $\sum_{i,j<N} (\tau_{i}\cdot D_{\tau_{j}}\eta)(\tau_{j}\cdot D_{\tau_{i}}\eta) = \sum_{i,j<N}\frac{\partial\eta^{i}}{\partial x_{j}}\frac{\partial\eta^{j}}{\partial x_{i}}$,\\
(xii) \begin{multline*}- \text{trace} ((\nabla \eta)^2) + \abs{\stackrel{\rightarrow}{n}\cdot\nabla\eta}^2 + 
2(\stackrel{\rightarrow}{n},\stackrel{\rightarrow}{n}\cdot(\nabla\eta)^2) - (\stackrel{\rightarrow}{n},\stackrel{\rightarrow}{n}\cdot\nabla\eta)^2\\ = 
 \abs{\stackrel{\rightarrow}{n}\cdot\nabla\eta}^2  - (\stackrel{\rightarrow}{n},\stackrel{\rightarrow}{n}\cdot\nabla\eta)^2 + 
2(\stackrel{\rightarrow}{n},\stackrel{\rightarrow}{n}\cdot(\nabla\eta)^2) - \text{trace} ((\nabla \eta)^2)
\\= \sum_{i=1}^{N-1}\abs{(D_{\tau_{i}}\eta)^{\perp}}^2 + (2\sum_{k}\frac{\partial\eta^{N}}{\partial x_{k}}\frac{\partial\eta^{k}}{\partial x_{N}} -\sum_{i,k}\frac{\partial\eta^{i}}{\partial x_{k}}\frac{\partial\eta^{k}}{\partial x_{i}} )\\
=\sum_{i=1}^{N-1}\abs{(D_{\tau_{i}}\eta)^{\perp}}^2  + (\abs{\frac{\partial\eta^{N}}{\partial x_{N}}}^2 -
\sum_{i,j<N}\frac{\partial\eta^{i}}{\partial x_{j}}\frac{\partial\eta^{j}}{\partial x_{i}})\\=\sum_{i=1}^{N-1}\abs{(D_{\tau_{i}}\eta)^{\perp}}^2  + 
\left((\stackrel{\rightarrow}{n},\stackrel{\rightarrow}{n}\cdot\nabla\eta)^2 - \sum_{i,j<N} (\tau_{i}\cdot D_{\tau_{j}}\eta)(\tau_{j}\cdot D_{\tau_{i}}\eta)\right).
\end{multline*}
Thus from (\ref{svreduced}), we find that
\begin{multline*}
\lim_{\varepsilon\rightarrow 0}\delta^{2}E_{\varepsilon}(u^{\varepsilon},\eta,\zeta) = 2\sigma\int_{\Gamma} \text{div}^{\Gamma}\zeta + (\text{div}^{\Gamma}\eta)^2 + 
\sum_{i=1}^{N-1}\abs{(D_{\tau_{i}}\eta)^{\perp}}^2  \\+ \left((\stackrel{\rightarrow}{n},\stackrel{\rightarrow}{n}\cdot\nabla\eta)^2 - \sum_{i,j<N} (\tau_{i}\cdot D_{\tau_{j}}\eta)(\tau_{j}\cdot D_{\tau_{i}}\eta)\right)\\
= \delta^{2}E(\Gamma,\eta,\zeta) + 2\sigma\int_{\Gamma} (\stackrel{\rightarrow}{n},\stackrel{\rightarrow}{n}\cdot\nabla\eta)^2.
\end{multline*}
The proof of our theorem is now complete.
\end{proof}
\begin{proof}[Proof of Theorem \ref{Morse}]
For any vector field $V$ defined on $\Gamma$ and is normal to $\Gamma$, we also denote by $V$ its extension to $\Omega$
in such a way that $(\stackrel{\rightarrow}{n}, \stackrel{\rightarrow}{n}\cdot \nabla V) =0.$ As a consequence, the second term on the right hand side of (\ref{discrep}) drops. Let $V$ and $W$ be vector fields normal to $\Gamma$. Then, let
\begin{equation*}
 \Phi_{V, t} = x + tV(x), \Phi_{W,t} = x + t W(x)
\end{equation*}
and
\begin{equation*}
 v_{\e} = u^{\e} (\Phi^{-1}_{V,t}), w_{\e} = u^{\e} (\Phi^{-1}_{W,t}).  
\end{equation*}
The proof is based on (\ref{discrep}) together with the following claims.
\begin{claim} (polarization)
\begin{equation}
B_{\e}(u^{\e}) (\partial_{t} v_{\e}(0), \partial_{t} w_{\e}(0)) = B(u) (V, W) + o(1). 
\label{polarident}
\end{equation}
\label{polari}
\end{claim}
\begin{claim} (injectivity)
The map $V\longmapsto \partial_{t} v_{\e}(0)$ is linear and one-to-one for $\e$ small. 
\label{indep}
\end{claim}
Now, having these claims, we can complete the proof of Theorem \ref{Morse}, following the arguments in the proof of Theorem 1.1 in \cite{Serfaty}. 
By definition of $n^{+}$, if $n^{+}$ is finite, we can find $n^{+}$ linearly independent vector fields $V^{1},\cdots, V^{n}$ which are defined on $\Gamma$ and normal
to $\Gamma$ such that
the quadratic function $Q(u)$ restricted to the space they span is positive, i.e., 
\begin{equation}
\min_{\sum_{i=1}^{n^{+}} a^2_{i}=1} Q(u) (\sum_{i=1}^{n^{+}} a_{i} V^{i})>0. 
\label{posi1}
\end{equation}
Denote $V^{i}_{\e} = \partial_{t} v^{i}_{\e} (0) = \left.\frac{d}{dt}\right\rvert_{t=0} u^{\e} \left(\left(x+ tV^{i}(x)\right)^{-1}\right)$. In view of (\ref{polarident}), we have for all $a_{i}$
\begin{equation*}
\lim_{\e\rightarrow 0}  Q_{\e}(u^{\e}) (\sum_{i=1}^{n^{+}} a_{i} V_{\e}^{i}) =  Q(u) (\sum_{i=1}^{n^{+}} a_{i} V^{i}) 
\end{equation*}
and the convergence is uniform with respect to $(a_{i})$ such that $\sum_{i=1}^{n^{+}} a^2_{i} =1$. Finally, we deduce from (\ref{posi1}) that for $\e$ small enough
\begin{equation}
\min_{\sum_{i=1}^{n^{+}} a^2_{i}=1} Q_{\e}(u^{\e}) (\sum_{i=1}^{n^{+}} a_{i} V_{\e}^{i})>0. 
\label{posi2}
\end{equation}
By Claim \ref{indep} and the linear independence of $V^{i}$, $V^{i}_{\e}$ are linearly independent for $\e$ small. Therefore, the $V^{i}_{\e}$ span a space of dimension $n^{+}$. This proves
that $D^{2} E_{\e}(u^{\e})$ has at least $n^{+}$ positive eigenvalues and thus $n^{+}_{\e}\geq n^{+}.$ Observe that if $n^{+} =+\infty$ then we can apply
the previous argument on subspaces of arbitrarily large finite dimension, and find that $n^{+}_{\e}$ is also $+\infty$ for $\e$ small. The same arguments work
for $n^{-}_{\e}$ and $n^{-}$.\\
\h {\it We now prove Claim \ref{polari}}. Indeed, using (\ref{discrep}), we see that
\begin{equation*}
 \left.\frac{d^2}{dt^{2}}\right\rvert_{t=0} E_{\e} (v_{\e}) = Q(u)(V)
\end{equation*}
or, equivalently
\begin{equation*}
 \lim_{\e\rightarrow 0} B_{\e} (u^{\e}) (\partial_{t} v_{\e}(0), \partial_{t}v_{\e}(0)) = B(u)(V, V).
\end{equation*}
Therefore
\begin{equation}
 B_{\e} (u^{\e}) (\partial_{t} v_{\e}(0), \partial_{t}v_{\e}(0)) = B(u)(V, V) + o(1).
\label{Beq}
\end{equation}
Applying (\ref{Beq}) to $V+W$ and $V-W$, we get
\begin{equation*}
 B_{\e} (u^{\e}) (\partial_{t} v_{\e}(0) + \partial_{t} w_{\e}(0), \partial_{t}v_{\e}(0) + \partial_{t} w_{\e}(0)) = B(u)(V + W, V + W) + o(1)
\end{equation*}
and
\begin{equation*}
 B_{\e} (u^{\e}) (\partial_{t} v_{\e}(0) - \partial_{t} w_{\e}(0), \partial_{t}v_{\e}(0) - \partial_{t} w_{\e}(0)) = B(u)(V - W, V - W) + o(1).
\end{equation*}
Subtracting these two relations, we obtain (\ref{polarident}).\\
\h {\it Finally, we prove Claim \ref{indep}}. Recall that $\Phi_{V,t}(x) = x + tV(x).$
Note that, for each $x$, we have 
\begin{equation*}
     x = \Phi_{V, t} (\Phi_{V,t}^{-1}(x)) = \Phi_{V,t}^{-1}(x) + t V (\Phi_{V,t}^{-1}(x)).
\end{equation*}
Hence 
\begin{equation*}
 0 = \frac{d}{dt}(\Phi_{V, t}^{-1}(x)) + t \nabla  V(\Phi_{V,t}^{-1}(x))\cdot\frac{d}{dt}(\Phi_{V, t}^{-1}(x)).
\end{equation*}
Evaluating the above equation at $t=0$ and noting that $\Phi_{V,0}^{-1} (x) =x$, one obtains
\begin{equation*}
\left. \frac{d}{dt}\right\rvert_{t=0}(\Phi_{V, t}^{-1}(x)) = -V(x).
\end{equation*}
It is now clear that
\begin{equation}
 V\longmapsto \partial_{t} v_{\e}(0): = \left.\frac{d}{dt}\right\rvert_{t=0} u^{\e} \left(\left(x+ tV(x)\right)^{-1}\right) = -\nabla u^{\e}\cdot V
\end{equation}
is a linear map. Let $V= f \stackrel{\rightarrow}{n}$ be a normal vector field to $\Gamma$. Suppose that $\nabla u^{\e}\cdot (f \stackrel{\rightarrow}{n}) =0$ for all $\e$ small. This implies that 
$\e\abs{\nabla u^{\e}\cdot \stackrel{\rightarrow}{n}}^2 f^2 =0$. Letting $\e\rightarrow 0$ and using (\ref{res}), we find that $2\sigma f^2 \equiv 0$ on $\Gamma$. Therefore $f=0$ and $V=0$.
This proves our claim.
\end{proof}
\begin{proof}[Proof of Theorem \ref{mMorse}]
Under the condition (C1) and the fact that $u^{\e}$ are critical points of $E_{\e}$, the work of Hutchinson and Tonegawa \cite[Theorem 1]{HT} (see also \cite{Tone}) showed that, in the sense of 
Radon measures, \begin{equation}\left(\frac{\varepsilon\abs{\nabla u^{\varepsilon}}^2}{2} +\frac{W(u^{\varepsilon})}{\varepsilon}\right)dx 
\rightharpoonup 2m\sigma d\mathcal{H}^{N-1}\lfloor\Gamma, \label{msingle} \end{equation} 
where $m$ is an integer-valued function defined on $\Gamma$. Furthermore, we have equipartition of energy, i.e.,
in the sense of Radon measures
\begin{equation}
 \abs{\frac{\varepsilon \abs{\nabla u^{\varepsilon}}^2}{2}-\frac{W(u^{\varepsilon})}{\varepsilon}} dx\rightharpoonup 0.
\label{vanishing}\end{equation}
Because $u^{\e}$ are critical points of $E_{\e}$, $\Gamma$ is a stationary varifold. Now, it follows from the connectivity of $\Gamma$ from (C2) and the
Constancy Theorem \cite[Theorem 41.1]{Simon} that $m$ must be a constant. From the constancy of $m$ and limiting equipartition 
of energy (\ref{vanishing}), it can be proved that (see, e.g., \cite[Equation (3.5)]{L})
\begin{equation}
 \varepsilon \nabla u^{\varepsilon}\otimes\nabla u^{\varepsilon} dx\rightharpoonup 2m\sigma \stackrel{\rightarrow}{n}
\otimes\stackrel{\rightarrow}{n} \mathcal{H}^{N-1}\lfloor\Gamma.\label{mres}
\end{equation}
Now, we proceed as in the proof of Theorem \ref{Morse} which used the result obtained in Theorem \ref{mainthm}. First, as in the proof of Theorem \ref{mainthm},
we have (\ref{svep}). Letting $\e\rightarrow 0$ in (\ref{svep}), using (\ref{res}), and computing as in the proof of Theorem \ref{mainthm}, we obtain 
 \begin{equation}\lim_{\varepsilon\rightarrow 0}\delta^{2} E_{\varepsilon}(u^{\varepsilon},\eta,\zeta) =
 m\delta^{2}E(\Gamma,\eta,\zeta) + 2m\sigma\int_{\Gamma} (\stackrel{\rightarrow}{n},\stackrel{\rightarrow}{n}\cdot\nabla\eta)^2\label{mdiscrep}\end{equation} 
for all smooth vector fields $\eta,\zeta\in (C_{c}^{1}(\Omega))^{N}$. \\
\h With (\ref{mdiscrep}), the proof of Theorem \ref{mMorse} can be completed similarly to that of Theorem \ref{Morse}.
\end{proof}
\begin{remark}
 There exist critical points $u^{\e}$ of $E_{\e}$ satisfying (\ref{msingle}) with any positive integer $m$ for minimal surfaces $\Gamma$ in bounded domains in
$\RR^{2}$ satisfying certain nondegeneracy conditions. See, e. g., \cite{dPKWToda}.
\end{remark}

{}


\begin{thebibliography}{xx} 
{\small
\bibitem{AC}Allen, S.; Cahn, J. W. A microscopic theory for antiphase boundary motion and its application to antiphase domain coarsening.
{\it Acta Metall.} {\bf 27} (1979) 1084--1095.
\bibitem{DFP} Dang, H.; Fife, P. C.; Peletier, L. A. Saddle solutions of the bistable diffusion equation.  
{\it Z. Angew. Math. Phys.}  {\bf 43}  (1992),  no. 6, 984--998.
\bibitem{dPKWToda}del Pino, M.; Kowalczyk, M.; Wei, J. The Toda system and clustering interfaces in the 
Allen-Cahn equation.  {\it Arch. Ration. Mech. Anal.}  {\bf 190}  (2008),  no. 1, 141--187.
\bibitem{dPKW} del Pino, M.; Kowalczyk M., Wei, J. C. Entire Solutions of the 
Allen-Cahn Equation and Complete Embedded Minimal Surfaces of Finite Total 
Curvature, preprint. 
\bibitem{Gui}Gui, C. Hamiltonian identities for elliptic partial differential equations.  {\it J. Funct. Anal.}  {\bf 254}  (2008),  no. 4, 904--933.
\bibitem{HT}Hutchinson, J. E.; Tonegawa, Y. Convergence of phase interfaces in the van der 
Waals-Cahn-Hilliard theory.  {\it Calc. Var. Partial Differential Equations}  {\bf 10}  (2000),  no. 1, 49--84.
\bibitem{L} Le, N. Q. A gamma-convergence approach to the Cahn-Hilliard equation.  {\it Calc. Var. Partial Differential Equations}  
{\bf 32}  (2008),  no. 4, 499--522.
\bibitem{LM} Luckhaus, S.; Modica, L. The Gibbs-Thompson relation within the gradient theory of phase transitions.  {\it Arch. Rational Mech. Anal.}  {\bf 107}  (1989),  no. 1, 71--83.
\bibitem{MM} Modica, L.; Mortola, S. Un esempio di $\Gamma \sp{-}$-convergenza. (Italian)  {\it Boll. Un. Mat. Ital. B (5)}  {\bf 14}  (1977), no. 1, 285--299. 
\bibitem{PR} Pacard, F.; Ritor\'{e}, M. From constant mean curvature hypersurfaces to the gradient theory of phase transitions. 
{\it J. Differential Geom.} {\bf 64} (2003), no. 3, 359--423. 
\bibitem{Res} Reshetnyak, Y. G. The weak convergence of completely additive vector functions on a 
set. {\it Siberian Math. J.} {\bf 9} (1968), 1039--1045; translated from Sibirskii Mathematicheskii Zhurnal {\bf 9} (1968), 1386-1394. 
\bibitem{RS}R\"{o}ger, M.; Sch\"{a}tzle, R. On a modified conjecture of De Giorgi.  {\it Math. Z.}  {\bf 254}  (2006),  no. 4, 675--714.
\bibitem{Serfaty} Serfaty, S. Stability in 2D Ginzburg-Landau passes to the limit.  {\it Indiana Univ. Math. J.}  {\bf 54}  (2005),  no. 1, 199--221
\bibitem{Simon} Simon, L. Lectures on geometric measure theory. Proceedings of the Centre for Mathematical Analysis, Australian National University, 3. Australian National University, Centre for Mathematical Analysis, Canberra, 1983.
\bibitem{Stern} Sternberg, P. The effect of a singular perturbation on nonconvex variational problems.  {\it Arch. Rational Mech. Anal.}  {\bf 101}  (1988),  no. 3, 209--260.
\bibitem{Tone} Tonegawa, Y. On stable critical points for a singular perturbation problem.  {\it Comm. Anal. Geom.}  {\bf 13}  (2005),  no. 2, 439--459
\bibitem{ToneDiff} Tonegawa, Y. A diffused interface whose chemical potential lies in a Sobolev space.  
{\it Ann. Sc. Norm. Super. Pisa Cl. Sci.} (5)  {\bf 4}  (2005),  no. 3, 487--510.
\bibitem{TW}Tonegawa, Y.; Wickramasekera, N. Stable phase interfaces in the van der Waals--Cahn--Hilliard theory, arXiv:1007.2060v1 [math.DG].
}
\end{thebibliography}
\end{document}